\documentclass[12pt]{amsart}
\usepackage{graphicx}
\usepackage{graphics}
\usepackage{amsmath}
\usepackage{amscd}
\usepackage{latexsym}

\begin{document}

\textwidth 6.2in
\textheight 7.6in
\evensidemargin .75in
\oddsidemargin.75in

\newtheorem{Thm}{Theorem}
\newtheorem{Lem}[Thm]{Lemma}
\newtheorem{Cor}[Thm]{Corollary}
\newtheorem{Prop}[Thm]{Proposition}
\newtheorem{Rm}{Remark}

\def\a{{\mathbb a}}
\def\C{{\mathbb C}}
\def\A{{\mathbb A}}
\def\B{{\mathbb B}}
\def\D{{\mathbb D}}
\def\E{{\mathbb E}}
\def\R{{\mathbb R}}
\def\P{{\mathbb P}}
\def\S{{\mathbb S}}
\def\Z{{\mathbb Z}}
\def\O{{\mathbb O}}
\def\H{{\mathbb H}}
\def\V{{\mathbb V}}
\def\Q{{\mathbb Q}}
\def\Cn{${\mathcal C}_n$}
\def\CM{\mathcal M}
\def\CG{\mathcal G}
\def\CH{\mathcal H}
\def\CT{\mathcal T}
\def\CF{\mathcal F}
\def\CA{\mathcal A}
\def\CB{\mathcal B}
\def\CD{\mathcal D}
\def\CP{\mathcal P}
\def\CS{\mathcal S}
\def\CZ{\mathcal Z}
\def\CE{\mathcal E}
\def\CL{\mathcal L}
\def\CV{\mathcal V}
\def\CW{\mathcal W}
\def\IC{\mathbb C}
\def\IF{\mathbb F}
\def\IK{\mathcal K}
\def\IL{\mathcal L}
\def\IP{\bf P}
\def\IR{\mathbb R}
\def\IZ{\mathbb Z}

\title{Nash homotopy spheres are standard}
\author{Selman Akbulut}
\thanks{The author is partially supported by NSF grant DMS 0905917}
\address{Department  of Mathematics, Michigan State University,  MI, 48824}
\email{akbulut@math.msu.edu }
\date{\today}
\begin{abstract} 
We prove that the infinite family of homotopy $4$-spheres constructed by Daniel Nash are all diffeomorphic to $S^4$.
\end{abstract}

\date{}
\maketitle

\setcounter{section}{-1}

\vspace{-.1in}

\section{Introduction}

In \cite{d}, D.Nash constructed an infinite family of smooth homotopy $4$-spheres  $\Sigma_{p,q,r,s}$ indexed by $(p,q,r,s)\in {\Bbb Z}^4$, and conjectured that they are possibly not diffeomorphic to  $S^4$. Here we prove that they are all diffeomorphic to $S^4$. In spirit, Nash's construction is an easier version of  the construction of the Akhmedov-Park in \cite{ap}, namely one starts with a standard manifold  $X^{4} =X_{1}\smile_{\partial} X_{2}$ which is a union of two basic pieces along their boundaries, then does the ``log transform" operations to some imbedded tories in both sides with the hope of getting an exotic copy of a known manifold $M^{4}$. In Nash's case $X$ is the double of $T_{0}^2\times T_{0}^2$ and $M$ is $S^4$,  in Akhmedov-Park case $X$ is the ``Cacime surface'' of \cite{ccm} (see \cite{a1} for decomposition of $X$) and $M$ is $S^2\times S^2$. 
 
\section{Log transform operation}

First we need to  recall the {\it  log transform} operation.  Let $X$ be a smooth $4$-manifold which contains a torus $T^{2}$ with the  trivial normal bundle $\nu (T^{2})\approx T^{2}\times B^2$.
  Let $\varphi_p$ $(p\geq 0)$ be the self-diffeomorphism of $T^3$ induced by the automorphism 
\begin{equation*}
\left(
\begin{array}{ccc}
1 &0 &0  \\
0 &p &-1  \\
0 &1 &0
\end{array}
\right)
\end{equation*}
of $H_1(S^1;\mathbf{Z})\oplus H_1(S^1;\mathbf{Z})\oplus H_1(S^1;\mathbf{Z})$ with the obvious basis $(a,b,c)$. 
 The operation of  removing $\nu (T) $ from $X$ and then regluing $T^2\times B^2 $ via $\varphi_{p} : S^1 \times T^2 \to \partial \nu (T)$ is called the  $p$  log-transform of $X$ along $T^2$.  In short we will refer this as $(a\times b,  b,  p)$ log transform. Figure 1 describes this as a handlebody operation (c.f. \cite{ay} and \cite{gs}).

 \begin{figure}[ht]  \begin{center}  
\includegraphics[width=.8\textwidth]{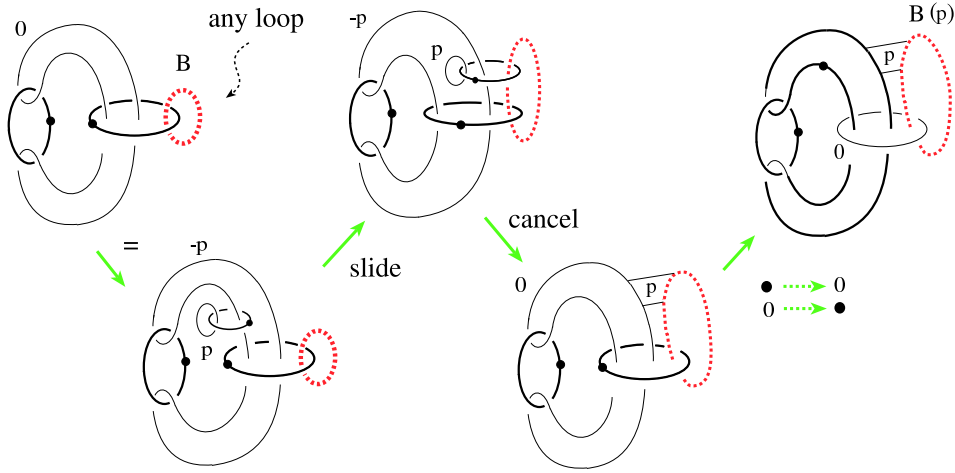}   
\caption{$p$ log-transform operation} 
\end{center}
\end{figure}

\section{Handlebody descriptions of  $T^2 \times T^2$ and $T^2_{0}\times T^{2}_{0}$}

We follow the recipe of \cite {a1} for drawing the surface bundles over surfaces. Here we have the simpler case of base and fiber are tori. Figure 2 describes the handlebody of $T^2$ and its thickening $T^2\times [0,1]$. Figure 3 is a handlebody of $T^2\times S^1$. Note that without the handle $v$ (the green curve) it is just $T^2_{0}\times S^1$. So  Figure 4 is the handlebody of $T^2\times T^2$. Note that without $2$-handles $u$ and $v$ (the green and blue curves)  Figure 4  it is just $T_{0}^2\times T_{0}^{2}$. By gradually converting the $1$-handles of Figure 4  from ``pair of balls" notation to the ``circle-with-dot'' notation of \cite{a2} we get the handlebody pictures of $T^2\times T^2$ and $T_{0}^2\times T_{0}^2$ in Figure 5 (the first and the second pictures). Then by an isotopy (indicated by the arrow)  we obtained the first picture of Figure 6 which is $T_{0}^2\times T_{0}^2$.

\section{Nash spheres are standard}
Let $X_{p,q}$ be the manifold obtained from $T_{0}^{2}\times T_{0}^{2}$ by  the  log-transformations $(a\times c,a,p)$ and 
 $(b\times c,b,q)$, where $a,b,c,d$ are the circle factors of $T_{0}^{2}\times T_{0}^{2}$ indicated of  Figure 6. Than Nash homotopy spheres are defined to be: 
$$\Sigma_{p,q,r,s}=X_{p,q}\smile_{\phi} -X_{r,s}$$
 where $\phi $ is the involution on  $\partial (T_{0}^{2}\times T_{0}^{2})$ flipping $T_{0}^{2}\times S^1$ and $S^{1}\times T^{2}_{0}$.  Notice that if $X^{r,s}$ is the manifold obtained from $T_{0}^{2}\times T_{0}^{2}$ by  the log transformations $(c\times a,c,r)$ and  $(d\times a,d,s)$ then we can identify:
 $$\Sigma_{p,q,r,s}=-X_{p,q}\smile_{id} X^{r,s}$$
 
 \begin{Thm} $\Sigma_{p,q,r,s}=S^4$
 \end{Thm}
 
 \proof  By using the description of the log-transform in Figure 1, we see that the second picture of Figure 6 is $X_{p,q}$. Now we will turn $X_{p,q}$ upside down and glue it to $X^{r,s}$ along its boundary. For this we take the image of the dual $2$-handle curves of $X_{p.q}$ (indicated by the dotted circles in Figure 6) by the diffeomorphism $\partial (X_{p,q})\stackrel{\approx}{\longrightarrow }\partial(X^{r,s})$,  and then attach $2$-handles to $X^{r,s}$ along the image of these curves. By reversing the log-transform process of Figure 6 we obtain the first picture of Figure 7, which is just $T^2_{0}\times T^2_{0}$, with the dual handle curves indicated! By an isotopy we obtain the second picture of Figure 7, also describing $T^2_{0}\times T^2_{0}$. Now again by using  the recipe of Figure 1 we  perform the $(c\times a,c,r)$ and  $(d\times a,d,s)$ log-transforms and obtain the first picture of Figure 8, which is $X^{r,s}$, with the dual $2$-handle curves of $X_{p,q}$ clearly visible in the picture. Now by attaching $2$-handles to top of $X^{r,s}$ along these curves we obtain $\Sigma_{p,q,r,s}$, which is also described by the first picture of Figure 8.  Here we should mention our convention: when a framing is not indicated in figures it means the zero framing. Now by the obvious handle slides and cancellations in Figure 8 $\leadsto$ Figure 9 we obtain $\sharp^{4} (S^2 \times B^2) $, and the four $3$-handles cancel these to give $S^4$ \qed.

\newpage

 \begin{figure}[ht]  \begin{center}  
\includegraphics[width=.6\textwidth]{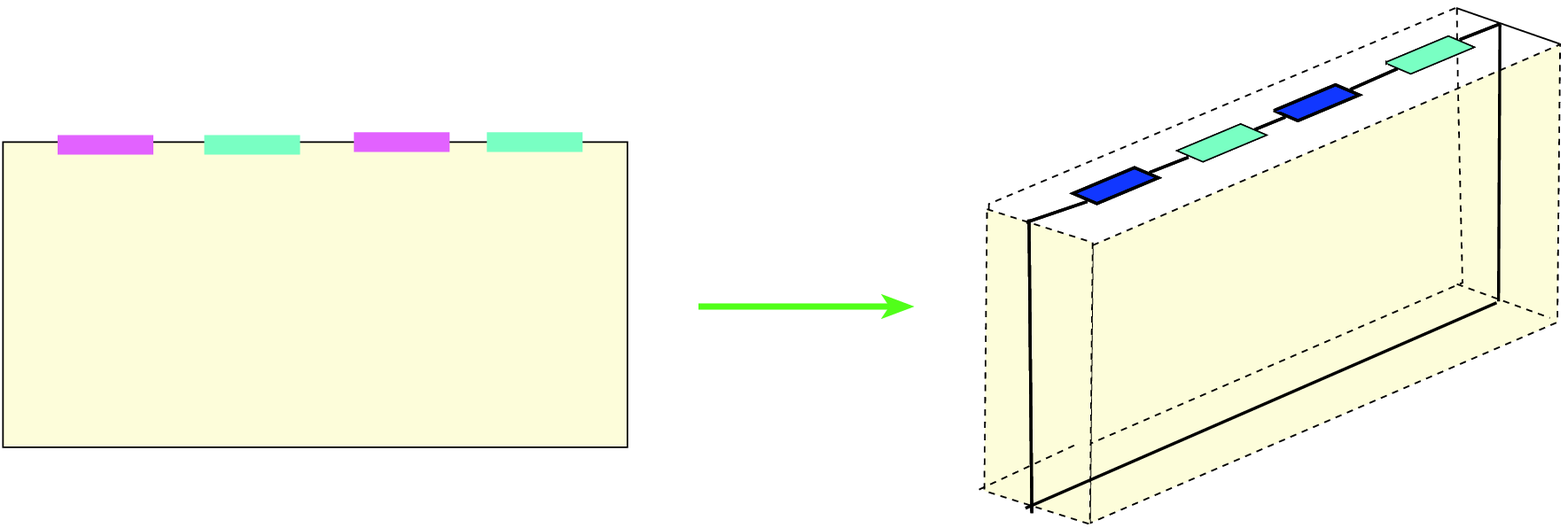}   
\caption{ $T^2\to T^2\times [0,1]$ } 
\end{center}
\end{figure}

 \begin{figure}[ht]  \begin{center}  
\includegraphics[width=.5\textwidth]{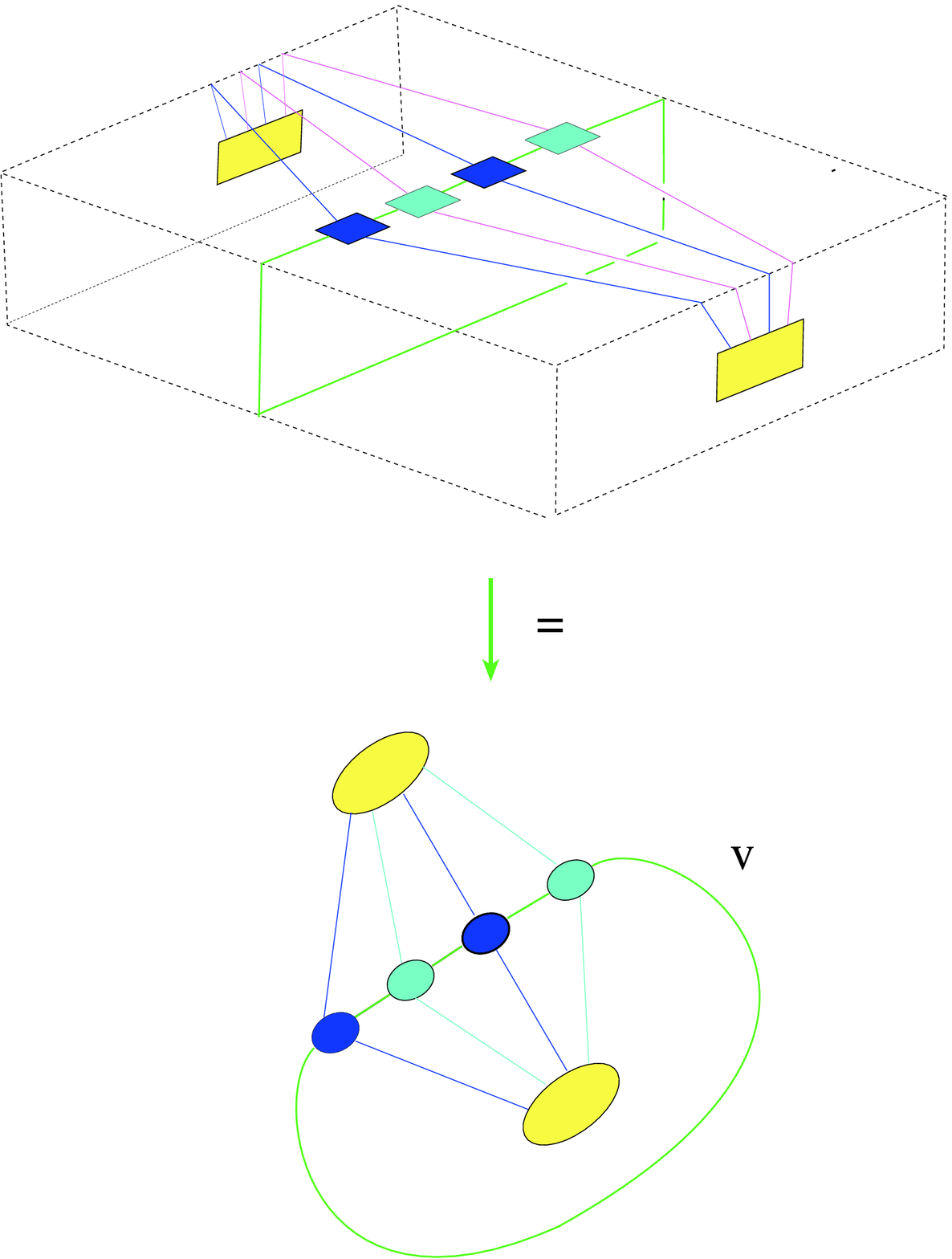}
\caption{$T^2\times S^1$}   
\end{center}
\end{figure}

 \begin{figure}[ht]  \begin{center}  
\includegraphics[width=.7\textwidth]{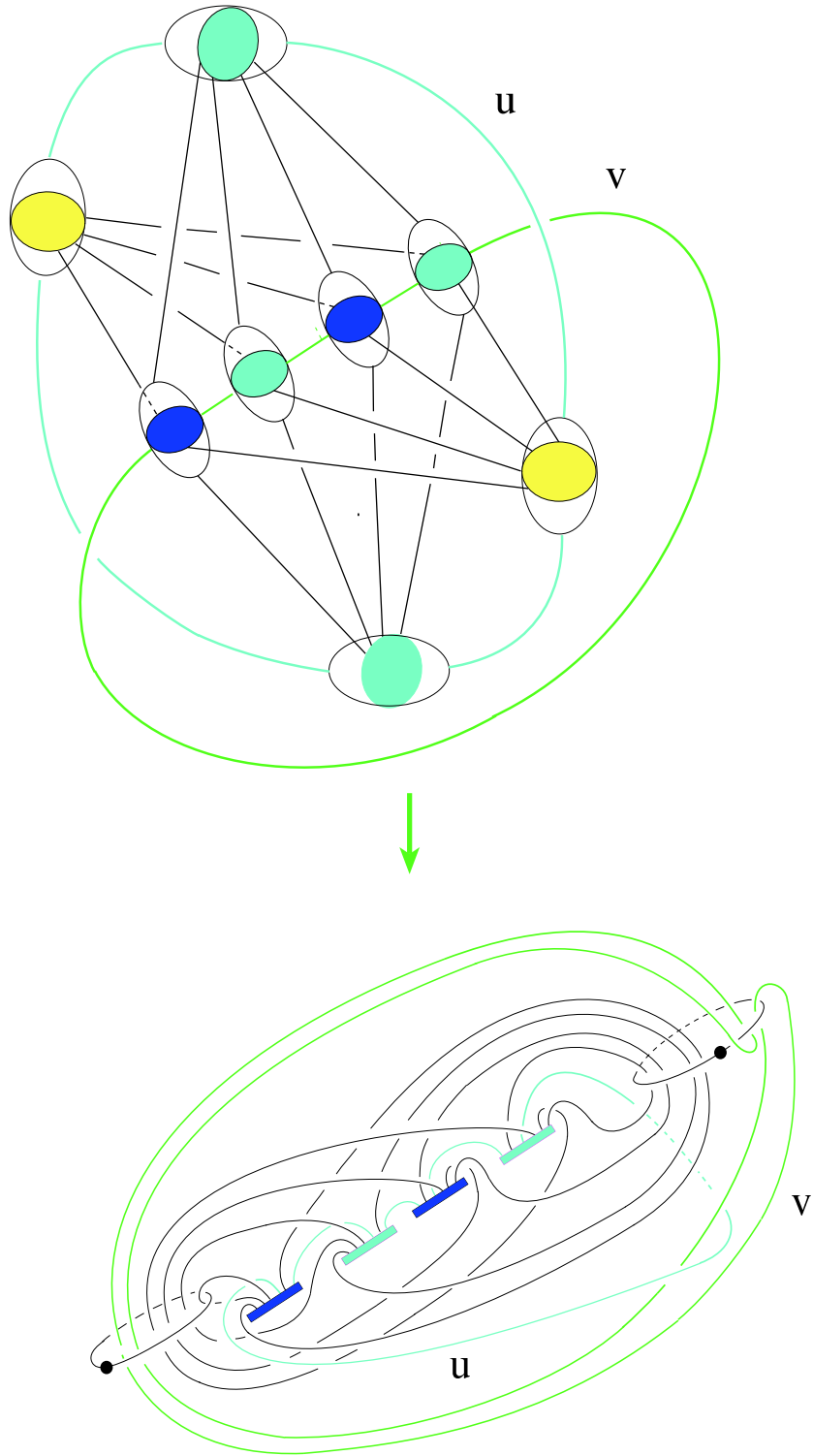}   
\caption{$T^2\times T^2$} 
\end{center}
\end{figure}

 \begin{figure}[ht]  \begin{center}  
\includegraphics[width=.8\textwidth]{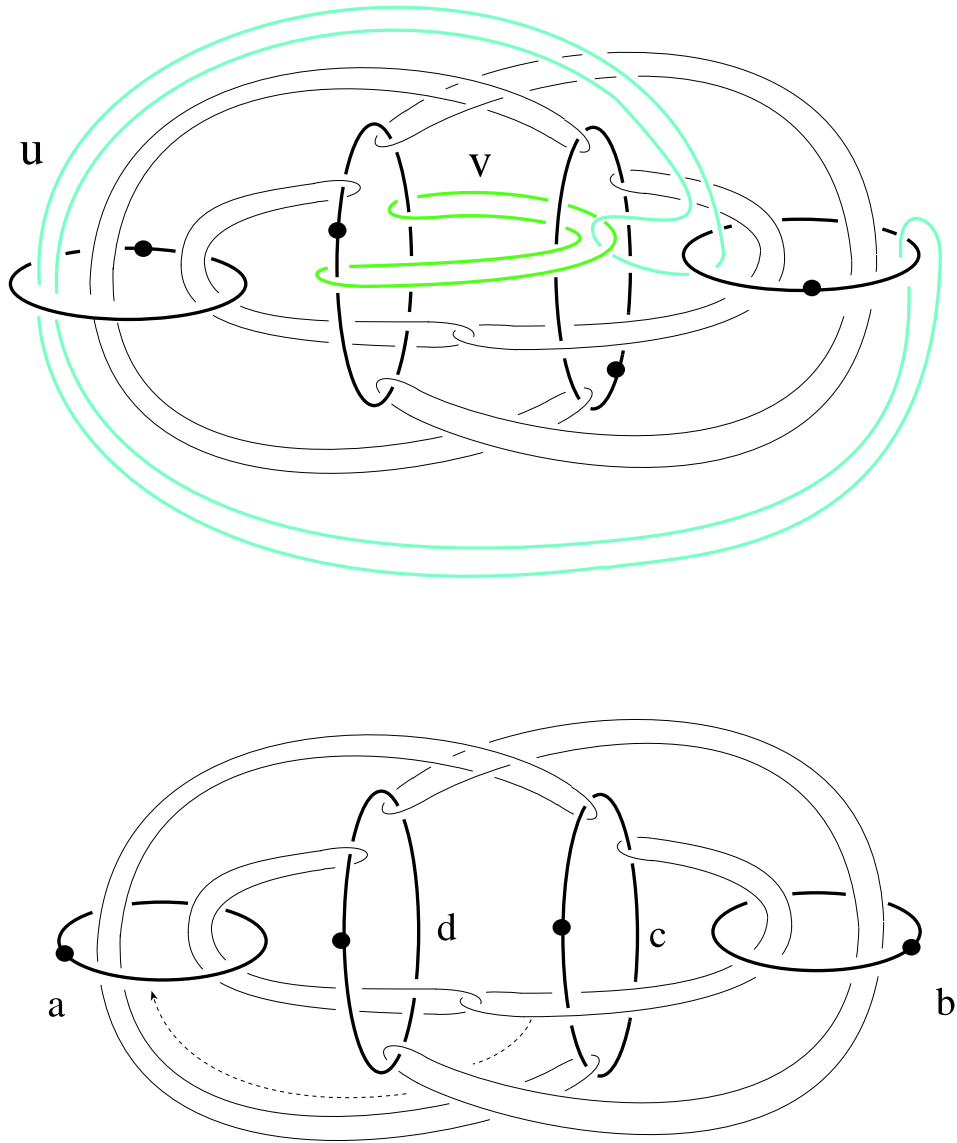}   
\caption{ $T^2\times T^2$ and $T^{2}_{0}\times T^{2}_{0}$} 
\end{center}
\end{figure} 

 \begin{figure}[ht]  \begin{center}  
\includegraphics[width=.8\textwidth]{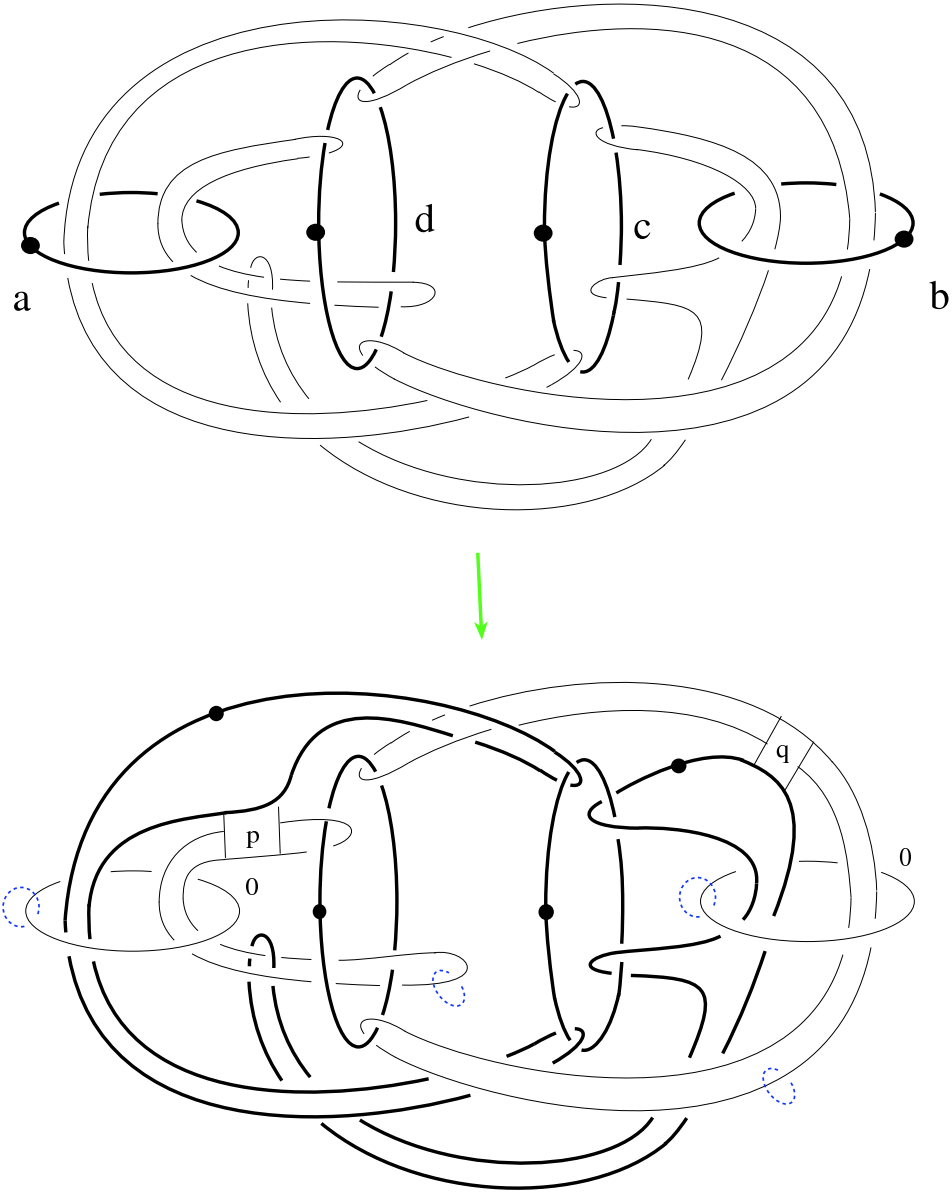}   
\caption{log-transformation $T_{0}^{2}\times T_{0}^{2} \to X_{p,q}$} 
\end{center}
\end{figure} 

\begin{figure}[ht]  \begin{center}  
\includegraphics[width=.8\textwidth]{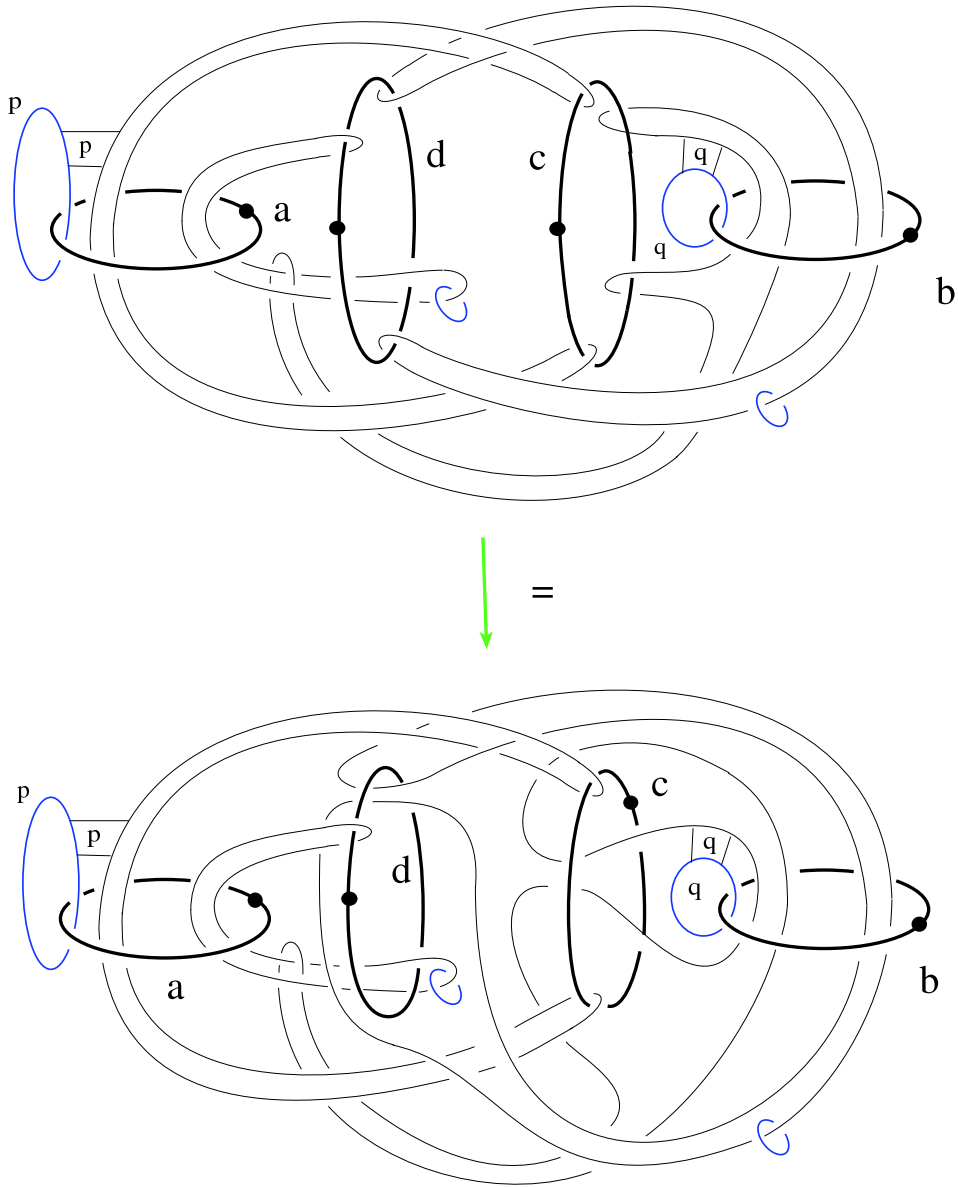}   
\caption{log-transformation undone} 
\end{center}
\end{figure} 

\begin{figure}[ht]  \begin{center}  
\includegraphics[width=.8\textwidth]{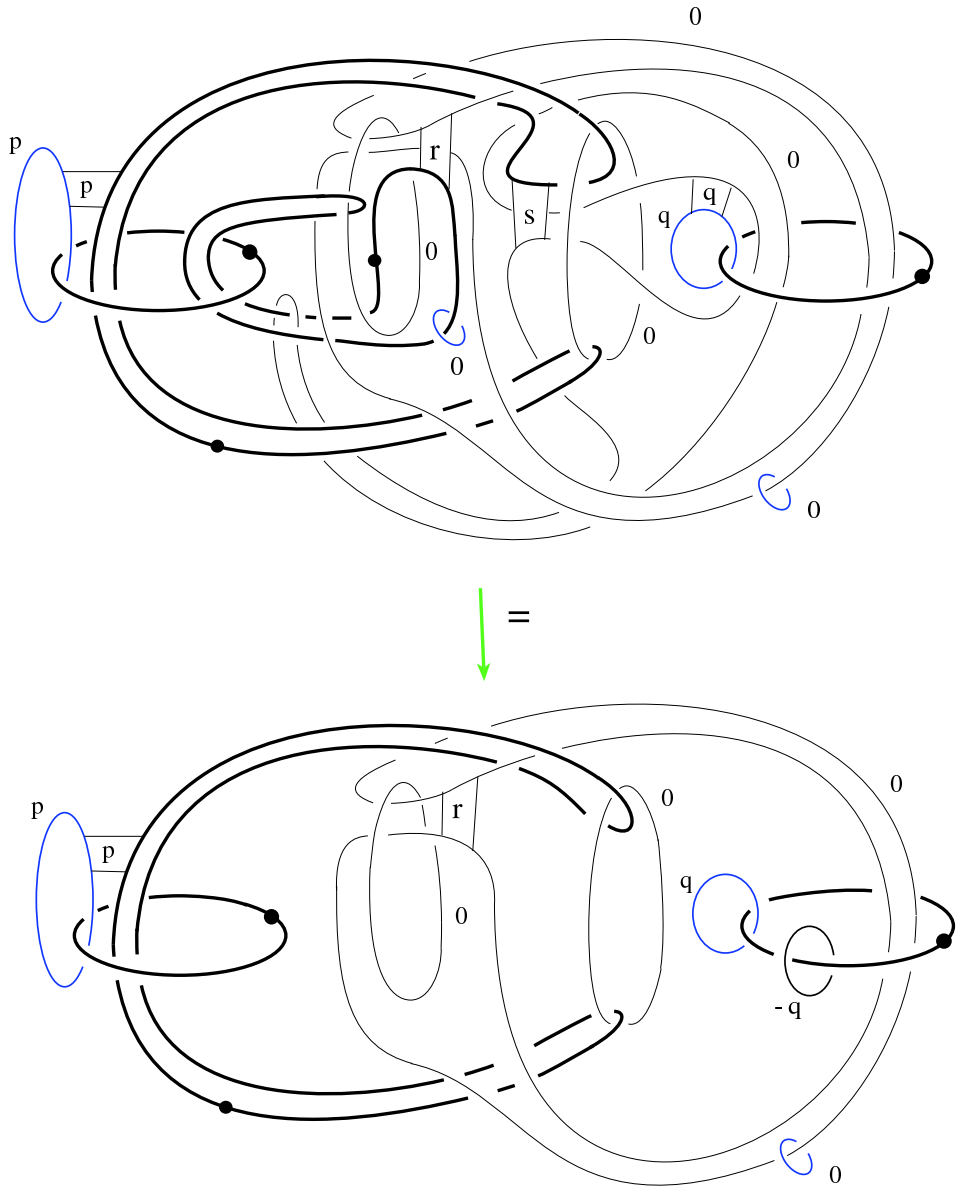}   
\caption{$\Sigma_{p,q,r,s}$} 
\end{center}
\end{figure} 

\begin{figure}[ht]  \begin{center}  
\includegraphics[width=.6\textwidth]{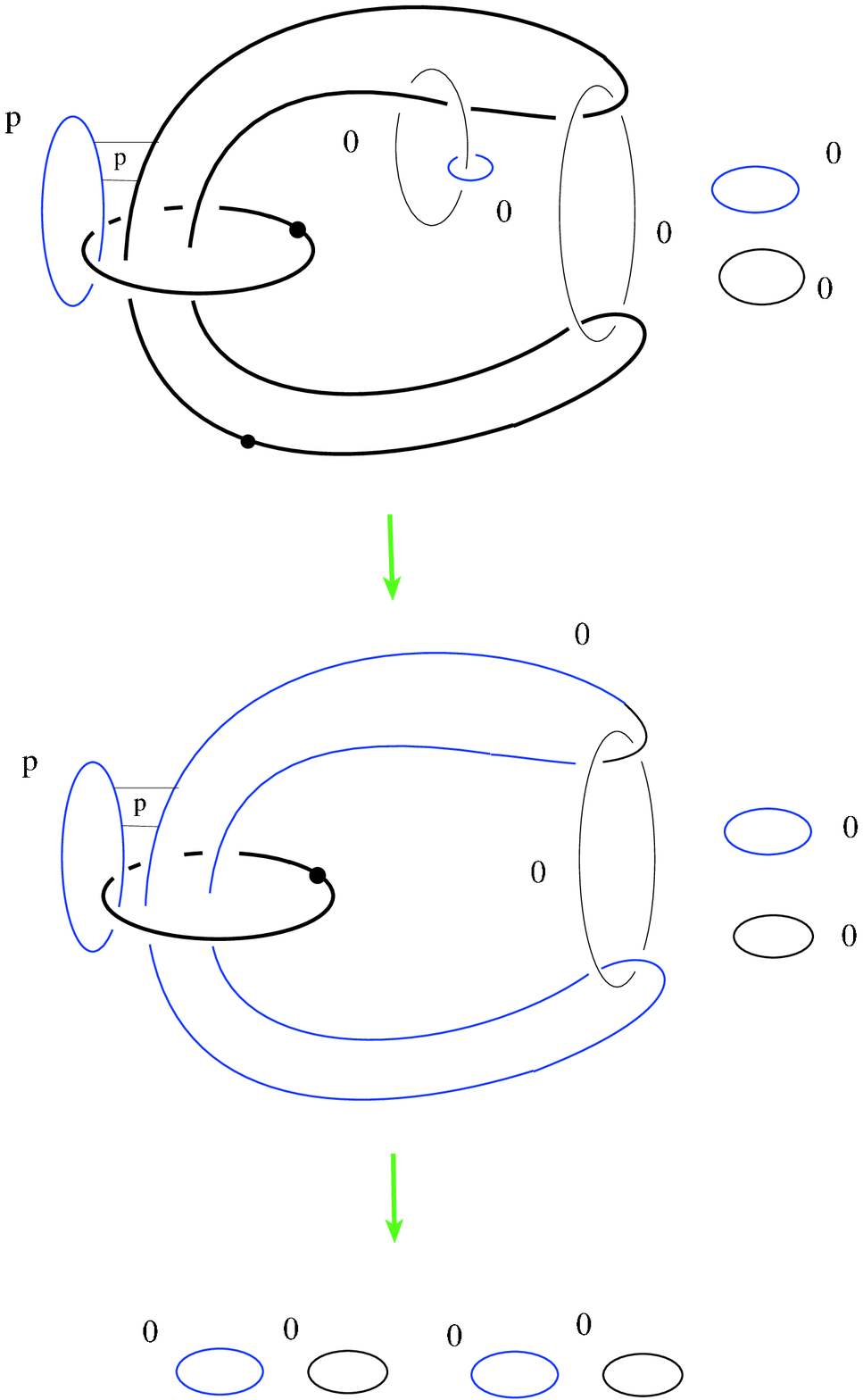}   
\caption{} 
\end{center}
\end{figure}

\end{document}